\documentclass{amsart}

\usepackage{latexsym}
\usepackage{epsfig}
\usepackage{amsthm}
\usepackage{amssymb}
\usepackage{amsmath}
\usepackage{amscd}
\usepackage{xcolor}
\usepackage{cleveref}

\newtheorem {Theorem}   {Theorem}

\newtheorem{lemma}[subsection]{Lemma}

\theoremstyle{remark}
\newtheorem{Remark}[Theorem]{Remark}

\newtheorem {Definition}[Theorem]{Definition}

\newtheorem{Conjecture}[Theorem]{Conjecture}

\newtheorem {proposition}[Theorem]{Proposition}
\newtheorem {corollary}[Theorem]{Corollary}
\expandafter\chardef\csname pre amssym.def at\endcsname=\the\catcode`\@
\catcode`\@=11
\def\undefine#1{\let#1\undefined}
\def\newsymbol#1#2#3#4#5{\let\next@\relax
 \ifnum#2=\@ne\let\next@\msafam@\else
 \ifnum#2=\tw@\let\next@\msbfam@\fi\fi
 \mathchardef#1="#3\next@#4#5}
\def\mathhexbox@#1#2#3{\relax
 \ifmmode\mathpalette{}{\m@th\mathchar"#1#2#3}%
 \else\leavevmode\hbox{$\m@th\mathchar"#1#2#3$}\fi}
\def\hexnumber@#1{\ifcase#1 0\or 1\or 2\or 3\or 4\or 5\or 6\or 7\or 8\or
 9\or A\or B\or C\or D\or E\or F\fi}

\font\teneufm=eufm10
\font\seveneufm=eufm7
\font\fiveeufm=eufm5
\newfam\eufmfam
\textfont\eufmfam=\teneufm
\scriptfont\eufmfam=\seveneufm
\scriptscriptfont\eufmfam=\fiveeufm

\catcode`\@=\csname pre amssym.def at\endcsname

\newcommand{\R}{{\bf R}}

\newcommand{\om}{{\omega}}

\newcommand{\pb}{\mathrm{pb}} 

\newcommand{\dd}{\mathrm{d}}

\usepackage{amssymb}
\usepackage[all]{xy}
\usepackage{enumerate}

\def    \R  {{\mathbb R}}

\def\qedbox{$\square$}%
\def\qed{\ifmmode\qedbox\else\unskip\ \hglue0mm\hfill
     \qedbox\smallskip\goodbreak\fi}%
\def\QED{\qed\goodbreak\vskip6pt}%

\begin{document}






\title[Smooth covers on symplectic manifolds]
{Smooth covers on symplectic manifolds}

\author[Fran\c{c}ois Lalonde]{Fran\c{c}ois Lalonde}
\author[Jordan Payette]{Jordan Payette}

\address{D\'epartement de math\'ematiques et de statistique, Universit\'e de Montr\'eal; 
D\'{e}partement de
math\'{e}matiques et de statistique, Universit\'{e} de Montr\'{e}al.}
\email{lalonde@dms.umontreal.ca; payettej@dms.umontreal.ca}

\thanks{The first author is supported by a Canada Research Chair, a NSERC grant OGP 0092913
(Canada) and a FQRNT grant ER-1199 (Qu\'ebec); the second author is supported by a Graham Bell fellowship from NSERC (Canada) }

\date{}

\begin{abstract}
   In this article, we first introduce the notion of a {\it continuous cover} of a manifold parametrised by any compact manifold $T$ endowed with a mass 1 volume-form. We prove that any such cover admits a partition of unity where the sum is replaced by integrals. When the cover is smooth, we then generalize Polterovich's notion of Poisson non-commutativity to such a context in order to get a more natural definition of non-commutativity and to be in a position where one can compare various invariants of symplectic manifolds. The main theorem of this article states that the discrete Poisson bracket invariant of Polterovich is equal to our smooth version of it, as it does not depend on the nature or dimension of the parameter space $T$. As a consequence, the Poisson-bracket invariant of a symplectic manifold can be computed either in the discrete category or in the smooth one, that is to say either by summing or integrating. The latter is in general more amenable to calculations, so that, in some sense, our result is in the spirit of the De Rham theorem equating simplicial cohomology and De Rham cohomology. We finally study the Poisson-bracket invariant associated to coverings by symplectic balls of capacity $c$, exhibiting some of its properties as the capacity $c$ varies. We end with some positive and negative speculations on the relation between uncertainty phase transitions and critical values of the Poisson bracket, which was the motivation behind this article.

\end{abstract}

\maketitle


\section{Introduction} 

      In mathematics, the notion of partition of unity is fundamental since it is the concept that distinguishes $C^{\infty}$ geometry from analytic geometry. In the first case, where partitions of unity apply, most objects can be decomposed in local parts, while in the second case where partitions of unity do not apply, most objects are intrinsically global and indecomposable.
      
      It is therefore of some importance to push that notion as far as we can in order to make it more natural and applicable. Our first observation is that the right context in which one should consider partitions of unity is in the continuous category (or possibly in the measurable category if one were able to make sense of that concept for families of open sets). So here continuous covers by open subsets of a given smooth manifold $M$ will be parametrised by any smooth compact manifold (possibly with boundary) $T$ endowed with a volume-form $dt$ of total mass $1$; for simplicity we shall refer to those pairs $(T, dt)$ as 'probability spaces'. We will first prove that for any compact manifold $M$, any such continuously parametrised cover of $M$ admits a smooth partition of unity made of smooth functions.
      
        Concentrating on an arbitrary symplectic manifold $(M, \om)$, the covers that we will consider will be made of smooth families of symplectically embedded balls of a given capacity $c= \pi r^2$ indexed by a measure space $(T, d\mu = \dd t)$. Here is our first theorem: the level of Poisson non-commutativity, as defined by Polterovich in the discrete case of partitions of unity, can be generalised to the case of our families of covers and associated partitions of unity; morever the number that we get in this general case, which depends {\it a priori} on the probability space, actually does not, being equal to the number associated to the corresponding discrete setting.

        Our second theorem is that if one considers the function $f: [0, c_{max}] \to [0, \infty]$  that assigns to $c$ the Polterovich's level of non-commutativity of the covers made of symplectically embedded balls of capacity $c$, as generalised by us in the smooth setting, then this function enjoys the following two properties:
        
        \medskip
        \noindent
        1) $f$ is non-increasing, and
        
         \medskip
        \noindent
        2) $f$ is upper semi-continuous and left-continuous.        
        \medskip

       We end this paper with a question concerning the relation, for a given symplectic manifold $(M, \omega)$, between critical values of the Poisson-bracket invariant as the capacity $c$ of the ball varies, and the critical values (or ``phase transition'') depending on $c$ of the topology of the infinite dimensional space of symplectic embeddings of the standard ball of capacity $c$ 
into $(M, \omega)$.  
       
       \vspace{6pt} \noindent {\it Acknowledgements}. Both authors are deeply grateful to Lev Buhovsky for suggesting and proving that the $T$-parameter spaces of smooth covers can be reduced to one-dimensional families. Although we present a different proof here, that now includes the discrete case, his idea has had a significant impact on the first version of this paper. The second author would like to thank Dominique Rathel-Fournier for inspiring discussions.

 \section{Continuous and smooth covers}
 
     Throughout this article, ``smooth" means ``of class $C^r$" for some arbitrary fixed $r \ge 1$ and $T$ is a compact smooth manifold of finite dimension endowed with a measure $\mu$ of total volume $1$ coming from a volume-form $\dd t$. The following definition is far more restrictive than the one that we have in mind, but it will be enough for the purpose of this article.
     
     \begin{Definition}
     Let $M$ be a closed smooth manifold of dimension $n$. Let $U$ be a bounded open subset of Euclidean space $\mathbb{R}^{n}$ whose boundary is smooth, so that the closure of $U$ admits an open neighbourhood smoothly diffeomorphic to $U$.
     A {\it continous cover} of $M$ of  type $(T,U)$ is a continuous map
     $$
     G: T \times U \to M $$
     such that
     
     \begin{enumerate}
     \item for each $t \in T$, the map $G_t$ is a smooth embedding of $U$ to $M$ that can be extended to a smooth embedding of some ({\it a priori} $t$-dependent) collar neighbourhood of $U$ (and therefore to the closed set $\bar{U}$), and
     
     \item the images of $U$ as $t$ runs over the parameter space $T$, cover $M$.
     \end{enumerate}
     
     \end{Definition} \vspace{6pt}
     
     Note that, in general, the topology of $U$ could change within the $T$-family. However, to simplify the presentation, we restrict ourselves to a fixed $U$ -- this is what we had in mind in the sentence preceding this definition. 
     A {\it smooth cover} is defined in the same way, but now requiring that $G$ be a smooth map.
               
       \begin{Definition} A {\it partition of unity} $F$ subordinated to a continuous cover $G$ is a smooth function
       $$
       \tilde{F}: T \times U \to [0, \infty)
       $$
       such that 
       
       \begin{enumerate}
       \item each $\tilde{F}_t : U \to \R$ is a smooth function with  (compact) support in $U$,
       
       \item the closure of the union $\bigcup_{t \in T} \, \mathrm{supp} (\tilde{F}_t)$ is contained in $U$, and
       
       \item for every $x \in M$,
       $$
       \int_T F_t(x) dt = 1,
       $$
        \end{enumerate}
       
\noindent
   where the smooth function $F_t : M \to \R$ is the pushforward of $\tilde{F}_t$ to $M$ using $G_t$, extended by zero outside the image of $G_t$;  in other words, it is $F_t(x) = \tilde{F}(t, G^{-1}_t(x))$.
       
       \noindent The notation $F < G$ expresses that $F$ is a partition of unity subordinated to the cover $G$.
       \end{Definition}

          \begin{Remark} Condition (2) plays a role in the proofs of a few results of this paper by allowing us to deform $U$ a little while keeping a given $F$ fixed; we were not able to come up with arguments working without this condition. Note that we recover the usual notion of partition of unity by taking $T$ to be a finite set of points with the counting measure.
                    \end{Remark} \vspace{0pt}

     \begin{Theorem}  Each continuous cover admits a partition of unity. 
      \end{Theorem}
     
     \proof
  Let $G$ be a continuous cover of $M$ of type $(T,U)$. The general idea of the proof is to replace $G$ by a finite open cover $G'$ of $M$, to consider a partition of unity subordinated to the latter and to use it to construct a partition of unity subordinated to $G$.
  
Cover $U$ by open balls such that their closure is always included inside $U$. Now, push forward this cover to $M$ using each $G_{t}$; the collection of all of these images as $t$ varies in $T$ forms an open cover of $M$ by sets diffeomorphic to the ball. Since $M$ is compact, there exists a finite subcover. Each open set in this subcover comes from some $G_t$, where $t$ is an element of a finite set $T' \subset T$. 

For each $t \in T'$, consider the (finite) collection $C_t$ of open balls inside $U$ whose image under $G_t$ belongs to the aforementioned subcover, so that the latter can be expressed as $G' := \{G_t(V) \, : \, t \in T', V \in C_t\}$. Since the closure of each ball $V \in C_t$ is contained in the open set $U$, by continuity of $G$ there is an open set $B_V \subset T$ centred at $t$ such that $G_t(\bar{V}) \subset G_{\tau}(U)$ for all $\tau \in B_V$, so that each $G_{\tau}^{-1} \circ G_t : \bar{V} \to U$ is defined and is a diffeomorphism onto its image. The intersection $B_t = \cap_{V \in C_t} B_V$ contains $t$ and is open since $C_t$ is finite. For each $t \in T'$, consider a smooth nonnegative bump function $\rho_t$ supported in $B_t$ whose integral over $T$ equals $1$.

There exists a smooth partition of unity $\Phi = \{ \phi_V : M \to [0,1] \, | \, V \in \cup_{t \in T'} C_t \}$ on $M$ subordinated to the finite open cover $G'$. For $t \in T'$ and $V \in C_t$, the real-valued function $\tilde{F}_V(\tau,u) := \rho_t(\tau) \phi_V(G(\tau,u))$ defined on $T \times U$ is supported in $B_t \times U$, is smooth in both $u$ and $\tau$ and satisfies $\int_T \tilde{F}_V(\tau, G_{\tau}^{-1}(x))d\tau = \phi_V(x)$. It easily follows that the function
\[ \tilde{F} : T \times U \to \mathbb{R} : (\tau, u) \mapsto \sum_{t \in T'} \sum_{V \in C_t} \tilde{F}_V(\tau, u) \]
is a partition of unity subordinated to $G$. 

\hfill $\square$ 
\bigskip

     \section{The $pb$ invariant and Poisson non-commutativity}
     
        Leonid Polterovich \cite{P1,P2,PR} introduced recently the notion of the level of Poisson non-com\-mutativity of a given classical (i.e finite) covering of a symplectic manifold. Here is the definition:
        
  \begin{Definition} Let $(M, \om)$ be a closed symplectic manifold and $\mathcal U$ a finite cover of $M$ by open subsets $U_1, \ldots, U_N$. For each partition of unity $F = (f_1, \ldots, f_N)$ subordinated to  $\mathcal U$, take the supremum of $\| \{\Sigma_i a_i f_i, \Sigma_j b_j f_j\} \|$ when the $N$-tuples of coefficients $(a_i)$ and $(b_i)$ run through the $N$-cube $[-1,1]^N$, where the bracket is the Poisson bracket and the norm is the $C^0$-supremum norm. Then take the infimum over all partitions of unity subordinated to $\mathcal U$. This is by definition the $pb$ \textit{invariant} of $\mathcal U$. To summarize:
  \[  \pb(\mathcal U)  := \underset{F < \mathcal U}{\mathrm{inf}} \; \underset{(a_i), (b_i) \in [-1,1]^N \, }{\mathrm{sup}} \, \left\| \left\{  \sum_i a_i f_i \, , \,  \sum_j b_j f_j    \right\} \right\| \; . \] \vspace{0pt}
  \end{Definition}

     Roughly speaking, this number is a measure of the least amount of ``symplectic interaction" that sets in a cover $\mathcal U$ can have. It is very plausible that such a number depends on the combinatorics of the cover, but also on the symplectic properties of the (intersections of the) open sets in the cover. To illustrate this point, observe that if $\mathcal U$ is an open cover made of only two open sets, then $\pb(\mathcal U) = 0$. A somewhat opposite result holds for covers constituted of displaceable open sets; let's recall that a subset $U$ of $M$ is {\it displaceable} if there is a Hamiltonian diffeomorphism $\phi$ such that $\phi(U) \cap U = \emptyset$. The main result of Polterovich in \cite{P1,PR} is that for such a cover, the number $\pb(\mathcal U)$ (multiplied by some finite number which measures the ``symplectic size" of the sets in the cover) is bounded from below by $(2N^2)^{-1}$. In particular, $\pb(\mathcal U) > 0$ in such a case. This result heavily relies on techniques in quantum and Floer homologies and in the theory of quasi-morphisms and quasi-states. Unfortunately, this lower bound depends on the cardinality $N$ of the open cover; as such, it does not show if one could use the $\pb$-invariant in order to assign to a given symplectic manifold a (strictly positive) number that might be interpreted as its level of Poisson non-commutativity. Nevertheless, Polterovich conjectured in \cite{P2} and \cite{PR} that for covers made of displaceable open sets, there should be a strictly positive lower bound for $pb$ independent of the cardinality of the cover, an extremely hard conjecture.
     
        One way of solving this problem might come from the extension of the $\pb$-invariant from finite covers to continuous or smooth covers. Indeed, such covers are morally limits of finite covers as the cardinality $N$ goes to infinity, so we can expect some relation between the minimal value of $\pb$ on such covers and the level of Poisson non-commutativity of the symplectic manifold.  This extension has the advantage that one may then compare the $\pb$ invariant for continuous/smooth covers to other quantities that also depend on continuous/smooth covers, such as the critical values at which families of symplectic balls undergo a ``phase transition''. We first need the following definition:
        
         \begin{Definition} Let $(M, \om)$ be a closed symplectic manifold and $G$ a continuous cover of $M$ of type $(T,U)$ by open subsets $G_t(U)$. For each partition of unity $F$ subordinated to $G$, take the supremum of $\| \{\int_N a(t) F(t) dt, \int_N b(t) F(t) dt\} \|$ over all coefficients (or \textit{weights}) $a$ and $b$ that are measurable functions defined on $T$ with $dt$-almost everywhere values in $[-1,1]$. Then take the infimum over all partitions of unity subordinated to $G$. This is by definition the $pb$ \textit{invariant} of $G$. To summarize:
         \[  \pb(G)  := \underset{F < G}{\mathrm{inf}} \; \underset{a, b : \, T \to [-1,1] \, \mbox{ \tiny{measurable}}}{\mathrm{sup}} \, \left\| \left\{  \int_T a(t) F_{t} dt \, , \,  \int_T b(t) F_{t} dt    \right\} \right\| \; . \] \vspace{0pt}
         
  \end{Definition}
  
     Note that we recover Polterovich's definition by replacing $T$ by a finite set of points. The following result shows that this pb-invariant is finite.
     
     \begin{lemma}
     Given a continuous cover $G$ of type $(T,U)$, there exists a partition of unity $F$ subordinated to $G$ whose pb-invariant $\pb \, F$ is finite\footnote{The pb-invariant $\pb \, F$ of a partition of unity $F$ is defined as above without the infimum over $F < G$. That is, $\pb \, G = \mathrm{inf}_{F < G} \, \pb \, F$.}.
     \end{lemma} \vspace{0pt}
     
     \proof 
     Consider the partition of unity $\tilde{F} : T \times U \to \mathbb{R}$ constructed in Theorem 4 above. Given a measurable function $a : T \to [-1,1]$, for any $x \in M$ we compute
     \[ \int_T a(\tau) F_{\tau}(x) d\tau = \sum_{t \in T'} \sum_{V \in C_t} \bar a_t \phi_V(x)  \, , \]
     where $\bar{a}_t := \int_T a(\tau) \rho_t(\tau)d\tau$ is a number whose value is in $[-1,1]$. It follows that
     \[ \left| \left\{ \int_T a(\tau) F_{\tau}(x) d\tau \, , \, \int_T b(\tau) F_{\tau}(x) d\tau \right\}(x)  \right| \le \sum_{V, W \in \cup_{t \in T'} C_t} |\{ \phi_V, \phi_W \}(x)| < \infty \]
     for any measurable functions $a, b : T \to [-1,1]$, which proves the claim. When the partition is smooth with respect to $t$ also, then the value is always finite. Indeed, since $F \in C^r(T \times M)$, it follows from Lebesgue's dominated convergence theorem that the function $a \cdot F := \int_T a(t) F_t dt \in C^r(M)$ is defined and also that $\{ a \cdot F, b \cdot F \} = \int_{T \times T} a(t) b(u) \{F_t, F_u \} dtdu < \infty$ when $a$ and $b$ are weights.

     \QED

 \noindent   We will be now working only with smooth covers. Therefore all partitions of unity $F$ satisfy $\pb \, F <  \infty.$

\medskip

          We recall a few facts taken from Polterovich's and Rosen's recent book \cite{PR}, since we will need them. Further informations are available in this book and in the references therein. The setting is the following (\cite{PR}, chapter 4):

\begin{itemize}
\item $(M^{2n}, \omega)$ a compact symplectic manifold;

\item $U \subset M$ an open set;

\item $H(U)$ the image of $\widetilde{Ham}(U)$ in $H := \widetilde{Ham}(M)$ under the map induced by the inclusion $U \subset M$;

\item $\phi$ : an element of $H(U)$;

\item $c$ : a (subadditive) spectral invariant on $H(U)$ (see the definition below);

\item $q(\phi) := c(\phi) + c(\phi^{-1})$, which is (almost) a norm on $H$;

\item $w(U) := \mathrm{sup}_{\phi \in H(U)} \, q(\phi)$ the spectral width of $U$ (which may be infinite).
\end{itemize} \vspace{6pt}

\begin{Definition}[\cite{PR}, 4.3.1]
A function $c : H \to \mathbb{R}$ is called a \textit{subadditive spectral invariant} if it satisfies the following axioms:
\begin{description}
\item[ Conjugation invariance ] $c(\phi \psi \phi^{-1}) = c(\psi) \; \forall \phi, \psi \in H$;

\item[ Subadditivity ] $c(\phi \psi) \le c(\phi) + c(\psi)$;

\item[ Stability ] 
\[  \int_{0}^1 \mathrm{min}(f_{t} - g_{t}) dt \le c(\phi) - c(\psi) \le \int_{0}^1 \mathrm{max}(f_{t} - g_{t}) dt \, ,  \] 

\noindent provided $\phi, \psi \in H$ are generated by normalized Hamiltonians $f$ and $g$, respectively;

\item[ Spectrality ] $c(\phi) \in \mathrm{spec}(\phi)$ for all nondegenerate elements $\phi \in H$.
\end{description}
\end{Definition} \vspace{0pt}

\begin{Remark} The first three properties of a spectral invariant are in practice the most important ones. However, from the spectrality axiom, one can show for instance that $w(U) < \infty$ whenever $U$ is displaceable; as such, the spectrality axiom is relevant in order to tie the spectral invariant with the symplectic topology of $M$. Let's mention that a spectral invariant exists on any closed symplectic manifold, as can be shown in the context of Hamiltonian Floer theory. \end{Remark}\vspace{0pt}

Given a Hamiltonian function $f \in C^{\infty}(M)$ generating the (autonomous) Hamiltonian diffeomorphism $\phi_{f} = \phi^1_{f}$ and a spectral invariant $c$, we can define the number
\[  \zeta(f) := \sigma(\phi_{f}) + \langle f \rangle \in \mathbb{R}   \]

\noindent where $\sigma(\phi_{f}) := \lim_{n \to \infty} \frac{1}{n}c(\phi^n_{f})$ (with $\sigma$ the {\it homogeneization} of $c$) and $\langle f \rangle := V^{-1} \, \int_{M} f \omega^n$ is the mean-value of $f$, where $V = \int_M \omega^n$ is the volume of the symplectic manifold $M$. The function $\zeta : C^{\infty}(M) \to \mathbb{R}$ is called the (\textit{partial symplectic}) {\it quasi-state} associated to $c$. It has some very important properties, among which:

\begin{description}
\item[ Normalization ] $\zeta(a) = a$ for any constant $a$;

\item[ Stability ] $\mathrm{min}_{M} (f-g) \le \zeta(f) - \zeta(g) \le \mathrm{max}_{M} (f-g)$;

\item[ Monotonicity ] If $f \ge g$ on $M$, then $\zeta(f) \ge \zeta(g)$;

\item[ Homogeneity ] If $s \in [0, \infty)$, then $\zeta(sf) = s \zeta(f)$;

\item[ Vanishing ] If the support of $f$ is displaceable, then $\zeta(f) = 0$ (this is a consequence of the spectrality axiom for $c$);

\item[ Quasi-subadditivity ] If $\{f, g \} = 0$, then $\zeta(f+g) \le \zeta(f) + \zeta(g)$.
\end{description} \vspace{12pt}

For $f, g \in C^{\infty}(M)$, define $S(f, g) = \mathrm{min} \{ w(\mathrm{supp} \, f) \; , \; w(\mathrm{supp} \, g) \} \in [0, \infty]$. It follows from Remark 8 that this number is finite whenever either $f$ or $g$ has displaceable support.\vspace{6pt}

\begin{Theorem}[\cite{EPZ}, 1.4 ; \cite{PR}, 4.6.1 ; the Poisson bracket inequality]
For every pair of functions $f, g \in C^{\infty}(M)$ such that $S(f,g) < \infty$,

\[  \Pi(f,g) := \left|  \zeta(f+g) - \zeta(f) - \zeta(g)  \right| \le \sqrt{2 S(f,g) \, \| \{ f, g \}  \| } \; . \] 
\end{Theorem} \vspace{12pt}

\noindent We see that $\Pi(f,g)$ measures the default of additivity of $\zeta$. In fact, this theorem implies:

\begin{description}
\item[ Partial quasi-linearity ] If $S(f,g) < \infty$ and if $\{f, g \} = 0$, then 
\[ \zeta(f+g) = \zeta(f) + \zeta(g) \; \mbox{ and } \zeta(s f) = s \zeta(f) \; \forall s \in \mathbb{R} \, . \]
\end{description} \vspace{12pt}

It is known that some symplectic manifolds admit a spectral invariant $c$ for which $S$ takes values in $[0, \infty)$, in which case $\zeta$ is a genuine symplectic quasi-state : it is a normalized, monotone and quasi-linear functional on the Poisson algebra $(C^{\infty}(M), \{-, - \})$. \vspace{6pt}

\begin{Theorem}[\cite{P1}, 3.1 ; \cite{PR}, 9.2.2]
Let $(M, \omega)$ be a symplectic manifold and consider a finite cover $U = \{ U_{1}, \dots, U_{N} \}$ of $M$ by displaceable open sets. Write $w(U) := \mathrm{max}_{i} \, w(U_{i}) < \infty$. Then
\[  \pb(U) \,  w(U) \; \ge \; \frac{1}{2N^2} \; .    \]
\end{Theorem} \vspace{6pt}

\textit{Proof} : \vspace{6pt}

Let $F$ be a partition of unity subordinated to $U$. Set
\[  G_{1} = F_{1}, \, G_{2} = F_{1} + F_{2}, \, \dots , \, G_{N} = F_{1} + \dots + F_{N} \, .   \]

\noindent Using Theorem 1 and the vanishing property of $\zeta$, one obtains the following estimate:
\begin{align}
\notag \left| \zeta(G_{k+1}) - \zeta(G_{k})  \right| &= \left| \zeta(G_{k} + F_{k+1} \, ) - \zeta(G_{k}) - \zeta(F_{k+1})   \right| \\
\notag & \le \sqrt{2 \, \mathrm{min} ( w(\mathrm{supp} \, G_{k}) \, , \, w(\mathrm{supp} \, F_{k+1})    )} \, \sqrt{\{ G_{k} , F_{k+1}   \}} \, . \end{align}

\noindent Using the definitions of $\pb (F)$ and of $w(U)$, one gets:
\begin{align}
\notag \left| \zeta(G_{k+1}) - \zeta(G_{k})  \right| &  \le \sqrt{2 \, w(U)} \, \sqrt{\pb (F) } \, .
\end{align}

\noindent This inequality holds for all $k$. Using the normalization and vanishing properties of $\zeta$ and applying the triangle inequality to a telescopic sum, one gets:
\begin{align}
\notag 1 & = \left| \zeta(1) - 0  \right| = \left| \zeta(G_{N}) - \zeta(G_{1})  \right| \le \sum_{k=1}^{N-1} \left| \zeta(G_{k+1}) - \zeta(G_{k})  \right| \\
\notag & \le \sum_{k=1}^{N-1}  \sqrt{2 \, w(U) \, \pb (F) } \le N \sqrt{2 \, w(U) \, \pb (F) } \, .
\end{align}

\noindent Since this is true for any $F < U$, the result easily follows.

\QED

A similar results holds in the context of smooth covers. We say that a smooth cover $G : T \times U \to (M, \omega)$ is made of displaceable sets if each set $G_t(\bar{U}) = \overline{G_t(U)} \subset (M, \omega)$ is displaceable (recall that we assume that $G_t$ extends as a smooth embedding to the closure $T \times \bar{U}$). In other words, not only is each $G_t(U)$ displaceable, but so is a small neighborhood of it too.

\begin{Theorem} For any smooth cover $G$ of type $(T,U)$ made of displaceable sets, there exists a constant $c = c(G) > 0$ such that
\[  \pb(G) \ge \; c(G) \; .    \]
\end{Theorem} \vspace{6pt}

\textit{Proof} : \vspace{6pt}

The proof morally consists in a coarse-graining of the smooth cover to a finite cover. Let $W_{1}, \dots, W_{N}$ be any exhaustion of the compact manifold $T$ by nested open sets with the following property: the sets $V_{1} = W_{1}$, $V_{2} = W_{2}-W_{1}$, ..., $V_{N} = W_{N} - W_{N-1}$ are such that for every $j$ the open set $U_j := \cup_{t \in V_{j}} \, \mathrm{Im}(G_{t})$ in $M$ is displaceable. Assume for the moment being that such sets $W_i$ exist. Notice that the sets $U_j$ cover $M$ and let $w(G) := \mathrm{sup}_j w(U_j) < \infty$. Now let $F$ be a partition of unity subordinated to $G$ and consider the functions $\int_{V_{1}} F_{t} dt$, ...\,, $\int_{V_{N}} F_{t} dt$ which form a partition of unity on $M$ subordinated to the $U_j$'s. As in the previous theorem, one estimates:

\begin{align}
\notag 1 &= \left| \zeta(1) - 0 \right| = \left| \zeta \left(\int_{W_{N}} F_{t}dt \right) - \zeta \left( \int_{W_{1}} F_{t}dt \right) \right| \\
\notag &\le \sum_{k=1}^{N-1} \left| \zeta \left(\int_{W_{k+1}} F_{t}dt \right) - \zeta \left( \int_{W_{k}} F_{t}dt \right) - \underset{0}{\underbrace{\zeta \left( \int_{V_{k+1}} F_{t}dt \right)}}\right| \\
\notag & \le \sum_{k=1}^{N-1} \sqrt{2 \, w(G) \, \pb(F)} \le N \sqrt{2 \, w(G) \, \pb(F)} \, .
\end{align} \vspace{0pt}

\noindent Since this is true for all $F < G$, and since $2N^2$ depends only on $G$ (through the choice of the $W_j$'s), the result follows with $c(G) := (2N^2 w(G))^{-1}$.

The sets $W_{j}$'s exist for the following reason. The closure of each $G_{t}(U)$ is a compact displaceable set, so that some open neighborhood $O_t$ of this set is displaceable. By the continuity of the cover $G$, for any $t$ there exists an open set $\{t \} \in Y_t \subset T$ such that $G(Y_t \times U) \subset O_t$. Since $T$ is compact, only a finite number of these $Y_t$ suffices to cover $T$, say $Y_1, \dots, Y_N$. Set $W_j = \cup_{k=1}^{j} Y_j$.  Since $V_j \subset Y_j$, the sets $G(V_j \times U)$ are indeed displaceable. This concludes the proof.

\QED 

It is natural to compare the $\pb$ invariant of different smooth covers of type $(T,U)$, especially if they are related to each other by a smooth family of smooth covers of the same type. This might help in understanding what is the 'optimal' way to cover a symplectic manifold $(M, \omega)$ by copies of a set $U$. We are led to the following definition which lies at the heart of this article:

\begin{Definition} A \textit{constraint} on smooth covers of $M$ of type $(T,U)$ is a set $C$ of such covers; the set of all smooth covers of type $(T,U)$ corresponds to the unconstrained case. Considering the $C^r$-Whitney topology on the space of smooth covers $G : T \times U \to M$, a \textit{constrained class of smooth covers of $M$ of type $(T,U)$} is defined as a connected component of the given constraint. We define the $\pb$ invariant of a (constrained) class $A$ as the infimum of $\pb(G)$ when $G$ runs over all smooth covers in $A$. \end{Definition}

As an instance of a constraint, we shall consider later on the one given by asking for each embedding $G_t : (U^{2n}, \omega_0) \hookrightarrow (M^{2n}, \omega)$ to be symplectic. The obvious difficulty with this last notion of pb invariant is that it intertwines four extrema: the supremum in the definition of the $C^0$-norm, the supremum over coefficients, the infimum over partitions of unity and the infimum over the smooth cover in the class. As a consequence of this difficulty, it is not clear if this number is strictly positive for every $M$, a problem which is related to Polterovich's conjecture; however, this number is now known to be positive for closed surfaces, as Polterovich's conjecture was recently proved valid in this context by Buhovsky, Tanny and Logunov \cite{BLT} and by the second author for genera $g \ge 1$ \cite{Pa}.

          \section{Equivalence of the smooth and discrete settings}

          This section is mainly devoted to the proof of Theorem 13 below which can be summarized as follows: the pb invariant of any class of $T$-covers is equal to the pb invariant of an affiliated class of discrete covers.\vspace{6pt}
          
          Fix a pair $(T,U)$. Any constraint $C$ of type $(T,U)$ determines the subset of \textit{constrained embeddings} 
          \[ C^* := \{ G_t : U  \hookrightarrow M \, | \, t \in T, \, G \in C \} \subset \mathrm{Emb}(U, M) \, .  \]
          
          \noindent Any section of the natural map $T \to \pi_0(T)$ -- which associates to $t \in T$ the connected component to which it belongs --  induces a well-defined, \textit{i.e.} section-independent, map $p_T : \pi_0(C) \to [\pi_0(T), C^*] \simeq \pi_0(C^*)^{\pi_0(T)}$. An element $A \in \pi_0(C)$ is just a constrained class of covers, and the element $A^* = p_T(A) \in \mathrm{Im}(p_T)$ corresponds to the $|\pi_0(T)|$ (not necessarily distinct) connected components of $C^*$ from which open sets the smooth covers in $A$ are built. Denote by $B$ the subset of $\pi_0(C^*)$ which is the image of $A^*$.
          Thus $B$ comprises sufficiently many open sets to cover the whole of $M$. \vspace{6pt}
          
          Let $\langle 1, n \rangle = [1,n] \cap \mathbb{N}$. Considering the natural map $q: C^* \to \pi_0(C^*)$, for $B \subset \pi_0(C^*)$ let $B' = q^{-1}(B) \subset C^*$. Assuming that $B'$ comprises enough open sets to cover $M$, define
          \[ \pb_{\mathrm{discrete}}(B) := \inf \, \{ \, \pb(G) \, | \, \exists n \in \mathbb{N}, \, G : \langle 1, n \rangle \to B' \mbox{ a cover of $M$ } \} \, . \]

     \noindent    To simplify the notations, we will, in the sequel, denote the set $B$ by the same symbol $A^*$.

          \begin{Theorem} [Equivalence smooth-discrete]
          Let $M$ be a symplectic manifold of dimension $2n$, $U$ an open subset of $\R^{2n}$ as mentioned above, and $T$ a compact manifold of strictly positive dimension endowed with a Lebesgue measure $\mu$ of total mass $1$. Consider a constraint $C$ on smooth covers of $M$ of type $(T,U)$, let $A \in \pi_0(C)$ be a constrained class of such covers and write $A^* = p_T(A) \subset \pi_0(C^*)$. Then
          \[ \pb(A) = \pb_{\mathrm{discrete}}(A^*) \; .   \]
          \end{Theorem}
          
          \proof
   We first prove $\pb(A) \ge \pb_{\mathrm{discrete}}(A^*)$. Let $G$ be a smooth cover of type $(T,U)$ in the constrained class $A$ and consider a smooth partition of unity $F < U$. By property (2) in the definition of a partition of unity and by continuity of $G$, we deduce that for each $t \in T$ there is an open set $t \in B_t \subset T$ such that $\mathrm{supp}(G_t^*(F_s)) \subset U$ for all $s \in B_t$. Since $T$ is compact, there is a finite set $T' = \{t_1, \dots, t_n\} \subset T$ such that the collection $B = \{B_{t_1}, \dots, B_{t_n}\}$ covers $T$. Consider a partition of unity $\rho = \{\rho_1, \dots, \rho_n\}$ on $T$ subordinated to $B$ and for each $t_i$ define
   \[ F'_i : M \to [0, \infty) : x \mapsto F_i(x) = \int_T \rho_i(t)F(x,t) dt \; . \]
   We observe that the collection $F' = \{F'_1, \dots, F'_n\}$ is a partition of unity on $M$ by smooth functions which is subordinated to the finite cover $G' := \left. G \right|_{T'}$ of $M$. We note that $\mathrm{Im}(G') \subset (A^*)'$, where we use a notation introduced just before the statement of the theorem. For $a' = \{a'_1, \dots, a'_n\} \subset [-1,1]$, the quantity $a := \sum_{i=1}^n a'_i \rho_i : T \to [-1,1]$ is a $T$-weight. For $a',b' \in [-1,1]^n$ we easily compute
   \begin{align}
   \notag  \left\{ \int_T a(t) F_t dt \, , \, \int_T b(u) F_u du \right\} &= \left\{ \sum_{i=1}^n a'_i F'_i \, , \, \sum_{j=1}^n b'_j F'_j \right\} \; .
   \end{align}
   Taking the suprema over weights thus yields $\pb(F) \ge \pb(F')$, while taking the infima over partitions of unities yields $\pb(G) \ge \pb(G')$. Taking the infima over covers in classes $A$ and $A^*$ finally yields $\pb(A) \ge \pb_{\mathrm{discrete}}(A^*)$. \vspace{6pt}
   
 We now prove $\pb(A) \le \pb_{\mathrm{discrete}}(A^*)$. Let $G' : \langle 1, n \rangle \to (A^*)'$ be a finite cover of $M$ and let $F' = \{F'_1, \dots, F'_n\}$ be a partition of unity subordinated to $G'$. Since $A^* = p_T(A)$, there exists a smooth cover $G''$ of $M$ of type $(T,U)$ in the constrained class $A \in \pi_0(C)$. Interpreting $A^* = \{A^*_1, \dots, A^*_m\}$ as a collection of connected components of $C^*$, for each connected component $A^*_i$ we can associate a point $t''_i \in T$ such that $G''_{t''_i} \in A^*_i$. From this association we can get an injective map $\langle 1, n \rangle \to T$ which associates to the integer $j$ a point $t'_j$ in the same connected component as the point $t''_i$, with $A^*_i \ni G'_j$. Call the image of this map $T' \subset T$.
 
 From these data we shall construct a smooth cover $G$ of type $(T,U)$ in the class $A$ which could act as a substitute for $G'$, in the sense that $\left. G \right|_{T'} = G'$. In fact, we shall define a smooth family $G_s$ of covers of type $(T,U)$ with $s \in [0,1]$ so that $G_0 = G''$ and $G_1 = G$, thereby illustrating that $G$ is indeed in the constrained class $A$. Fix a Riemannian metric on $T$. Observe that smoothly deforming $G''$ within $A$ if necessary, we can assume that $G''$ is constant in an $\epsilon$-neighbordhood of $T'$. If some connected component of $T$ contains none of the points $t'_j$, just set $G_s = G''$ on that component. For any other connected component of $T$, say the one containing $t''_i$, pick a Riemannian metric on it and consider disjoint embedded closed geodesic $\epsilon$-balls centred at the points $t'_j$. Outside the reunion of these balls, set again $G_s = G''$, whereas on the ball containing $t'_j$ define $G_s$ as follows. First choose a smooth path $g_j : [0, \epsilon] \to A^*_i$ such that $g_j(0) = G'_j$and $g_j(\epsilon) = G''(t'_j)$. Also pick a smooth function $\chi : [0, \epsilon] \to [0,1]$ such that $\chi(u) = 1$ if $u < \epsilon/3$ and $\chi(u) = 0$ is $u > 2\epsilon/3$. Denoting $r(p)$ the radial distance in the $j$-th ball of a point $p$ from $t'_j$, set on that ball $G_s(p) = g_j([1 - (1- \chi(s \epsilon)) \chi(r)]\epsilon)$.  This completely defines the family $G_s$ in the way we desired.
 
 We observe that $G$ is constant on an $(\epsilon/3)$-neighbourhood of each $t'_j$. For each $j$, pick a smooth positive function $\rho_j$ with support in the $(\epsilon/3)$-ball about $t'_j$ and which integrates to $1$. We define the smooth function $F : T \times M \to [0, \infty)$ as $F(t,m) = \sum_{j=1}^n \rho_j(t) F'_j(x)$. We easily verify that this is a smooth partition of unity subordinated to $G$.
 
 For any $T$-weight $a : T \to [-1,1]$, define $a' = (a'_1, \dots, a'_n) \in [-1,1]^n$ via $a'_j = \int_T a(t) \rho_j(t) dt$. For $T$-weights $a$ and $b$ we then easily compute
   \begin{align}
   \notag    \left\{ \sum_{i=1}^n a'_i F'_i \, , \, \sum_{j=1}^n b'_j F'_j \right\} &= \left\{ \int_T a(t) F_t dt \, , \, \int_T b(u) F_u du \right\} \; .
   \end{align}
   Taking the suprema over weights thus yields $\pb(F') \ge \pb(F)$, while taking the infima over partitions of unities yields $\pb(G') \ge \pb(G)$. Taking the infima over covers in classes $A^*$ and $A$ finally yields $\pb_{\mathrm{discrete}}(A^*) \ge \pb(A)$.

   \QED      \vspace{12pt}

   \section{Independence on the probability space}   
   
   The equivalence of the smooth and of the discrete settings suggests that the pb invariants might be independent from the underlying probability space $T$ parametrising the smooth covers. The purpose of this section is make this idea precise.
   
   \begin{proposition}
   Let $M$ be a symplectic manifold of dimension $2n$, $U$ an open subset of $\R^{2n}$ as mentioned above, and $T_1$ and $T_2$ be compact manifold of strictly positive dimension each endowed with a smooth volume form of total mass $1$. Consider constraints $C_1$ and $C_2$ on smooth covers of $M$ of type $(T_1,U)$ and $(T_2, u)$, respectively.  Let $A_i \in \pi_0(C_i)$, $i=1,2$, be constrained classes and assume that the corresponding sets of embeddings $(A_i^*)' \subset C^*_i \subset \mathrm{Emb}(U, M)$ coincide in the latter space. Then
          \[ \pb(A_1) = \pb(A_2) \; .   \]
   \end{proposition}
   
   \proof
   It follows from Theorem 13 that $\pb(A_i) = \pb_{\mathrm{discrete}}(A_i^*)$, $i=1,2$. Looking at the definition, $\pb_{\mathrm{discrete}}(A_i^*)$ only depends on the set $(A_i^*)'$, which is itself assumed to be independent from $i$.
   \QED

     Next we discuss special sorts of constraints which not only frequently appear in practice, but also for which the hypothesis in the previous proposition follows from a somewhat less stringent assumption.
        
        \begin{Definition}
        A constraint $C$ on covers of type $(T,U)$ is \textit{prime} if there exists a set $C' \subset \mathrm{Emb}(U, M)$ such that $G \in C$ if and only if $G_t \in C'$ for every $t \in T$. In other words, $C$ is prime if it is the largest constraint such that $C^* \subset C'$ (equivalently, $C^* = C'$).
        \end{Definition}\vspace{6pt}
        
        \noindent We point out that $C'$ thus admits sufficiently many open sets to cover the whole of $M$. Conversely, given a set $C' \subset \mathrm{Emb}(U, M)$ which admits sufficiently many open sets to cover $M$ and a probability space $T$, it is not guaranteed that there exists a constraint $C$ of covers of type $(T,U)$ (let alone a prime one) such that $C^* = C'$; this happens if $|\pi_0(C')| > |\pi_0(T)|$ and if no reunion of $|\pi_0(T)|$ connected components of $C'$ has sufficiently many open sets to cover $M$. In comparison, as long as $|\pi_0(C')|$ is finite, we can always find a discrete cover of $M$ made of open sets in $C'$. Note that this is however the only obstacle: given a set $C' \subset \mathrm{Emb}(U, M)$ such that there exists a smooth cover $G$ of $M$ of type $(T,U)$ with $G_* : \pi_0(T) \to \pi_0(C')$ well-defined and surjective, then $C' = C^*$ for some (prime) constraint $C$ on covers of type $(T,U)$.
        
      \begin{Definition}
        A prime constraint $C$ on covers of type $(T,U)$ with $C^* = C' \subset \mathrm{Emb}(U,M)$ is \textit{filled} if there is $G \in C$ such that the map $G_* : \pi_0(T) \to \pi_0(C')$ is surjective. By extension, we say that $C'$ is \textit{filled by $T$} if the associated prime constraint $C$ of type $(T,U)$ is filled.
        \end{Definition}\vspace{0pt}

     \begin{corollary}
   Let $M$ be a symplectic manifold of dimension $2n$, $U$ an open subset of $\R^{2n}$ as mentioned above, and $T_1$ and $T_2$ be compact manifold of strictly positive dimension each endowed with a smooth volume form of total mass $1$. Consider constraints $C_1$ and $C_2$ on smooth covers of $M$ of type $(T_1,U)$ and $(T_2, u)$, respectively.  Let $C' \subset \mathrm{Emb}(U,M)$ be filled by both $T_1$ and $T_2$ and consider the corresponding prime constraints $C_1$ and $C_2$. Let $A_i \in \pi_0(C_i)$, $i=1,2$, be constrained classes and assume that the corresponding sets of embeddings $(A_i^*)' \subset C^*_i \subset \mathrm{Emb}(U, M)$ coincide in the latter space. Then
          \[ \pb(A_1) = \pb(A_2) \; .   \]
     \end{corollary}
     
     \begin{Remark} \label{independence}
     For one application of this corollary, note that $C'$ is filled by any probability space $T$ of strictly positive dimension whenever $C'$ is connected and contains sufficiently many open sets to cover $M$. In that case $(A^*)'=C'$ for any $A \in \pi_0(C)$ (where $C$ is the prime and filled constraint associated with $C'$), since in fact $|\pi_0(C)|=1$. As a consequence, when $C'$ is not necessarily connected but each of its components contains sufficiently many embeddings to cover $M$, then the restriction of $\pb$ to prime constrained classes of covers parametrised by \textit{connected} $T$ comes from a function on $\pi_0(C')$.
     \end{Remark}

     \vspace{6pt}     
    \section{The behaviour of $\pb$ on symplectic balls}

     \bigskip \bigskip
       For the rest of this article, we only\footnote{The results of this section can however be easily adapted for star-shaped domain $U \subset \mathbb{R}^{2n}$.} consider $U = U(c) = B^{2n}(c)$, that is the standard symplectic ball capacity $c = \pi r^2$ (where $r$ is the radius). We also only consider (smooth) \textit{symplectic} covers, that is covers $G$ of type $(T,U)$ satisfying the symplectic prime constraint $C$ given as follows: $G \in C$ if $G_t \in C' = \mathrm{Emb}_{\omega}(U, M)$ for every $t \in T$. We shall write $U(c)$, $C(c)$ and $C'(c)$ when we want to stress the dependence on $c$. \vspace{6pt}
       
       Of special interest are the cases when $T = S^n$ for some $n \ge 1$. A constrained class $A$ of $C$ determines a connected component\footnote{It is still a conjecture, that we shall dub the \textit{symplectic camel conjecture}, whether $C'(c)$ is connected (whenever nonempty) when $(M, \omega)$ is compact and for any $c$.} $A' = p_T(A) \subset C'$, and determines in fact an element of the $n$-th homotopy group $\pi_n(A')$. Conversely, since $M$ is compact and using the fact that the group $\mathrm{Symp}(M, \omega)$ is $k$-transitive for all $k \in \mathbb{N}$, any element in $\pi_n(C')$ can be represented by some class $A \in \pi_0(C)$. The pb-invariants of symplectically constrained classes hence allow to probe the homotopic properties of $C'(c)$, properties which might change with $c$. Consequently, it appears important to better understand how the pb-invariants depend on the capacity $c$. This behavior of $\pb$ on $c$ is the main question raised in this paper.\vspace{6pt}
       
However, invoking  \Cref{independence} and again the $k$-transitivity of $\mathrm{Symp}(M, \omega)$, we deduce that for any connected probability space $T$ there is a bijective correspondence between $\pi_0(C)$ and $\pi_0(C')$. We can thus interpret the $\pb$ functional on smooth covers of type $(T,U(c))$ parametrised by connected spaces $T$ simply as a map $\pb : \pi_0(C'(c)) \to  [0, \infty)$, the latter being clearly independent from $T$. It therefore appears that the $\pb$-invariants can only probe the homotopy type of $C'(c)$ in a crude way. \vspace{6pt}
         
         Let $c_{max}$ denote the largest capacity a symplectic (open) ball embedded in $M$ can have; that can be much smaller than the one implied by the volume constraint $\mathrm{Vol}(U(c))  \le \mathrm{Vol}(M, \omega)$, according to the Non-Squeezing Theorem. For $0 < c < c' < c_{max}$ the obvious inclusion $U(c') \subset U(c)$ induces the restriction map $C'(c) \to C'(c')$ and hence also $r_{c,c'} : \pi_0(C'(c)) \to \pi_0(C'(c'))$.

          \begin{Definition} The \textit{tree of path-connected classes of symplectic embeddings of $U$ in $M$} is the set
          \[ \Psi(U,M) := \bigsqcup_{c \in (0, c_{max})} \{c\} \times \pi_0(C'(c)) \; .  \]
          A \textit{(short) branch of $\Psi(U,M)$} is a continuous path $\beta : (0, c_{\beta}) \to \Psi(U,M) : c \mapsto (c, A^*_{\beta}(c))$ such that $r_{c,c'}(A^*_{\beta}(c)) = A^*_{\beta}(c')$.
          \end{Definition} \vspace{0pt}
          
          We can therefore define a function $\pb : \Psi(U,M) \to [0, \infty)$ in the obvious way. Given a branch $\beta$ with domain $(0, c_{\beta})$, we can define a map $\pb_{\beta} = \pb \circ \beta : (0, c_{\beta}) \to [0, \infty)$.

      \begin{Theorem} Given any branch $\beta$, the function $\pb_{\beta}$ is non-increasing, upper semi-continuous and left-continuous.
      \end{Theorem}

 \proof 
 
  (a) Let us first show that the function is non-increasing. Fix $0 < c' < c < c_{\beta}$ and let $\epsilon > 0$. From the work done above and with the interpretation of $A^*_{\beta}(c)$ as a connected component of $C'(c)$, there exists a discrete cover $G' : \langle 1, n \rangle \to A^*_{\beta}(c')$ of $M$ such that $\pb(G') < \pb_{\beta}(c') + \epsilon$. We claim that this cover refines a cover $G : \langle 1, n \rangle \to A^*_{\beta}(c)$ of $M$; assuming this for the moment, we would then have
 \[ \pb_{\beta}(c) \le \pb(G) \le \pb(G') < \pb_{\beta}(c') + \epsilon \; . \] 
 As this holds for any $\epsilon > 0$, we get $\pb_{\beta}(c) \le \pb_{\beta}(c')$ \textit{i.e.} $\pb_{\beta}$ is non-increasing.
 
  To prove the existence of $G$, consider a symplectic embedding $B \in A^*_{\beta}(c)$. Since $\beta$ is a branch, the restriction $B'$ of $B$ to $U(c')$ is an embedding in $A^*_{\beta}(c')$; the latter space being a connected component of $C'(c')$ with respect to the Whitney $C^r$-topology, for each $j \in \langle 1, n \rangle$ there is smooth path of symplectic embeddings of $U(c')$ into $M$ joining $B'$ to $G'(j)$. By the symplectic isotopy extension theorem, each of these paths extends to a global symplectic isotopy on $M$, which thus sends $B$ to an embedding $G(j)$ of $U(c)$ into $M$. Clearly $G$ is a discrete cover of $M$ refined by $G'$.
     
   (b)  Now let us show that for every $c \in (0, c_{beta})$, the function $\pb_{\beta}$ is upper semi-continuous at $c$, \textit{i.e.} $\limsup_{c' \to c} \pb_{\beta}(c') \le \pb_{\beta}(c)$.
   
    On the one hand, it follows from part (a) that $\pb_{\beta}(c)$ is greater or equal to all limits of $f$ from the right. On the other hand, for any $\epsilon > 0$, there are a discrete cover $G$ representing $A^*_{\beta}(c)$ and a partition of unity $F < G$ such that $\pb(F) < \pb_{\beta}(c) + \epsilon$.  In fact, by our definition of a partition of unity, there is a strictly smaller capacity $c' < c$ such that the support of $F$ is compact inside the open ball $U(c') \subset U(c)$. Transporting the data to the restriction of the pair $(G,F)$ to $U(c'')$ for any $c'' \in [c',c]$, one gets 
   $$\pb_{\beta}(c'') \le \pb_{\beta}(c) + \epsilon.$$
   
  \noindent Since the choice of $c'$ indirectly depends on $\epsilon > 0$ through $F$, and might get as close to $c$ when $\epsilon$ approaches to zero, we do not get $\pb_{\beta}(c'') \le \pb_{\beta}(c)$ but only that $\pb_{\beta}$ is upper semi-continuous from the left. 
  
 (c) We wish to prove that $\pb_{\beta}$ is in fact left-continuous, that is to say that $\pb_{\beta}(c)$  is equal to the limit of $\pb_{\beta}(c')$ as $c'$ tends to $c$ from the left. Consider a sequence of capacities $c_i < c$ converging to $c$ with highest value $\lim \pb_{\beta}(c_i)$ (the value $\infty$ is not excluded).    This limit cannot be smaller than $\pb_{\beta}(c)$ because otherwise it would contradict the non-increasing property. However, by upper semi-continuity, it cannot be greater than $\pb_{\beta}(c)$. Therefore, it has to be equal to $\pb_{\beta}(c)$.
   
   \QED
   
   With regard to the continuity of the function $\pb_{\beta}$ associated to a branch $\beta$ it is not possible to be much more specific than the above Theorem, at least not when $\mathrm{dim} \, M = 2$. Indeed, in that case $c_{max} = \mathrm{Area}(M, \omega)$ and Moser's argument allows to prove that the space $C'(c) = \mathrm{Emb}_{\omega}(B^{2}(c), M)$ is connected whenever non-empty, so that there is only one maximal branch $\beta$. Polterovich's conjecture has recently been established in dimension two \cite{BLT}: in fact there is a universal constant $\gamma > 0$ such that $\pb_{\beta}(c)c > \gamma$ whenever $c \le c_{max}/2$. Using the invariance of the quantity $\pb_{\beta}(c)c$ upon pullback of the data under any symplectic covering map, this inequality holds even for $c > c_{max}/2$ when $M$ has genus $g \ge 1$ (\textit{c.f.} \cite{Pa}). However for $M = S^2$, by enlarging two opposite hemispheres one gets $\pb_{\beta}(c) = 0$ when $c > c_{max}/2$. Consequently $\pb_{\beta}$ is discontinuous on $S^2$, yet might be continuous on higher genus surfaces.

  \section{Phase transitions and the $\pb$ function}
  
      We conclude this paper by a few speculations since they disclose the main motivation behind this article. 
  
   The first ``phase transition'' discovered in Symplectic Topology is the following one:
          
     \begin{Theorem} (Anjos-Lalonde-Pinsonnault) In any ruled symplectic 4-manifold $(M, \om)$, there is a unique value $c_{crit}$ such that the infinite dimensional space $Emb(c,\om)$ of all symplectic embeddings of the standard closed ball of capacity $c$ in $M$ undergoes the following striking property: below $c_{crit}$, the space $Emb(c,\om)$ is homotopy equivalent to a finite dimensional manifold, while above that value, $Emb(c,\om)$ does not retract onto any finite dimensional manifold (or CW-complex) since it possesses non-trivial homology groups in dimension as high as one wishes. Below and above that critical value, the homotopy type stays the same.
     \end{Theorem}
     
      \medskip
   \noindent The reason for the term {\it phase transition} is still debatable, but there are several physical reasons, from Thermodynamics, to adopt that terminology. \vspace{6pt}

    \begin{Definition}  Given a closed symplectic manifold $(M, \omega)$, let us call an {\it uncertainty phase transition} any critical value $c$ at which the space of symplectic embeddings of balls of capacity $c$ into $(M, \omega)$ undergoes a change of its homotopy type.     \end{Definition}
    
    This terminology reflects the fact that a symplectically embedded ball quantifies the uncertaintly in the position and momentum of a (collection of) particles(s).
    
       The proof of the above theorem is quite indirect: one identifies all homology classes of symplectically embedded balls through the action of two groups on them: the full group of symplectic diffeomorphisms and the subset of these that preserve a given standard ball, the latter being viewed as the group of symplectic diffeomorphisms on the blow-up. Each of theses groups is computed by their action on a stratification of all compatible almost complex structures that realise holomophically some homology classes (essentially the homology classes that cut out the symplectic manifold in simple parts). Everything boils down to the behaviour of some J-curves in the given symplectic manifold for each $J$, generic or not (the non-generic ones playing the fundamental role since only the first stratum is generic). So, for instance, some homology class of symplectically embedded balls may disappear at some capacity $c_{crit}$ because the homology class of symplectic diffeomorphisms that preseve some standard ball of capacity $c$ in $M$ vanishes when $c$ crosses $c_{crit}$. It is conceivable that the class that vanished was supporting a covering that minimized the $pb$ at that level of capacity. We know that the dimension of the homology class, i.e. the dimension of the parametrizing space $T$, plays no role by our theorem on smooth-discrete equivalence. However, it is possible, that such a class, discretized or not, contained the optimal configuration of balls for a covering in order to mimimize $pb$.  Therefore the main question that drove us to study the $pb$ invariant in the smooth setting is:
       
\noindent
{\bf Question (Poisson-Uncertainty).} Is there a relation between the critical values of the Poisson bracket and the critical values (or phase transitions) of $Emb_{\omega}(B(c), M)$ as $c$ varies ?

   This is a natural question since the latter probes the topological changes in configurations of balls, while the former looks for $pb$-optimal configurations. We do not have in mind any direct sketch proving that there is a relation. So we must simply for the moment look at the facts. We have little material to work on, since the $pb$ conjecture has been proved (very recently) only for real surfaces, while the study of the topology of balls is known only in dimension $2$ and $4$ for ruled symplectic $4$-manifolds. Thus we may just examine the case of surfaces. In this case, there is no critical value for the phase transition, but there are for the pb-invariant, showing that the answer to the above question is negative in dimension $2$.

Small displaceable balls should not see the symplectic form, actually the space of (unparamatrised) symplectic balls below the uncertainty critical value retracts to the topology of the manifold itself for ruled symplectic 4-manifolds. This refines the symplectic camel conjecture for small capacities and it leads us to state the following conjecture:

\begin{Conjecture} (The Topology conjecture). The limit of the function $c \pb(c)$, as $c$ tends to zero, is a finite number, and depends only on the differential topology of the symplectic manifold.
\end{Conjecture} \vspace{6pt}

  Now, while the the Poisson-uncertainty question might have a positive answer in high dimensions,  we show here that the Poisson-uncertainty question has a negative answer in dimension $2$ for the sphere.

    To see this, let us consider the simple situation of $(M, \omega)$ being $S^2$ with its standard symplectic form, say of area $A$. As $M$ is a surface, it satisfies the symplectic camel conjecture, which is to say that the space $\mathrm{Emb}(c, \omega)$ is connected. The Poisson bracket function is then defined for any  $c \in (0, A)$. There exists on any closed symplectic manifold a spectral invariant $c$ such that $c(Id) = 0$, see Theorem 4.7.1 in \cite{PR}. It follows from that and the other properties of $c$ that the spectral width $w(U)$ of any subset $U \subset M$ satisfies $w(U) \le 4 e_H(U)$ where $e_H(U)$ is the Hofer displacement energy of $U$. For open sets in $S^2$, $e_H(U) = \mathrm{Area}(U)$ if this area is smaller than $A/2$ and $e_H(U) = \infty$ otherwise. In this context, Polterovich's conjecture (now a theorem on surfaces \cite{BLT,Pa}) states that there is a constant $C > 0$ such that for any, continuous or discrete, cover $G$ of $S^2$ by displaceable open sets, the inequality
\[ \pb(G)w(G) \ge C  \; \mbox{ holds }. \]

\noindent This implies that $\pb(c) e_H(U(c)) \ge C$. Thus when $c < A/2$, we have $\pb(c) \ge 2C/A$. However, we observe that $\pb(c) = 0$ whenever $c > A/2$: two symplectic balls of capacity $c > A/2$ suffice to cover $S^2$ and  the $\pb$-invariant of such a cover vanishes. Polterovich's conjecture hence goes against any claim that the Poisson bracket function $\pb(c)$ only has discontinuities when $\mathrm{Emb}(c, \omega)$ undergoes a transition in its homotopy type. 
\vspace{6pt}

As a concluding remark, we point out that our borrowings in the thermodynamical and statistical mechanical terminology are explained by our insight that tools from these subjects might play a role in the understanding of the symplectic problems we considered in this paper. The space of symplectically embedded balls can be understood as an infinite dimensional (pre)symplectic manifolds which is some sort of limit of finite dimensional ones. In this paper, continuous covers have also been understood as limits of discretes ones. It is a recurrent theme in statistical mechanics that systems with a very large number of degrees of freedom tend to behave in  universal and somewhat simpler ways.

 \end{document}